# A COMMUTATIVE DIAGRAM AMONG DISCRETE AND CONTINUOUS NEUMANN BOUNDARY OPTIMAL CONTROL PROBLEMS

**Domingo A. Tarzia**

CONICET - Depto. Matemática,
FCE, Univ. Austral, Paraguay 1950,
S2000FZF Rosario, Argentina.
Tel.: +54-341-5223093; Fax: +54-341-5223001
E-mail: DTarzia@austral.edu.ar

## ABSTRACT

We consider a bounded domain $\Omega \subset \mathbb{R}^n$ whose regular boundary $\Gamma = \partial\Omega = \Gamma_1 \cup \Gamma_2$ consists of the union of two disjoint portions $\Gamma_1$ and $\Gamma_2$ with positive measures. The convergence of a family of continuous Neumann boundary mixed elliptic optimal control problems ($P_\alpha$), governed by elliptic variational equalities, when the parameter $\alpha$ of the family (the heat transfer coefficient on the portion of the boundary $\Gamma_1$) goes to infinity was studied in Gariboldi - Tarzia, Adv. Diff. Eq. Control Processes, 1 (2008), 113-132, being the control variable the heat flux on the boundary $\Gamma_2$. It has been proved that the optimal control, and their corresponding system and adjoint system states are strongly convergent, in adequate functional spaces, to the optimal control, and the system and adjoint states of another Neumann boundary mixed elliptic optimal control problem ($P$) governed also by an elliptic variational equality with a different boundary condition on the portion of the boundary $\Gamma_1$.

We consider the discrete approximations ($P_{h\alpha}$) and ($P_h$) of the optimal control problems ($P_\alpha$) and ($P$) respectively, for each $h > 0$ and for each $\alpha > 0$, through the finite element method with Lagrange's triangles of type 1 with parameter $h$ (the longest side of the triangles). We also discrete the elliptic variational equalities which define the system and their adjoint system states, and the corresponding cost functional of the Neumann boundary optimal control problems ($P_\alpha$) and ($P$). The goal of this paper is to study the convergence of this family of discrete Neumann boundary mixed elliptic optimal control problems ($P_{h\alpha}$) when the parameter $\alpha$ goes to infinity. We prove the convergence of the discrete optimal controls, the discrete system and adjoint system states of the family ($P_{h\alpha}$) to the corresponding to the discrete Neumann boundary mixed elliptic optimal control problem ($P_h$) when $\alpha \to \infty$, for each $h > 0$, in adequate functional spaces. We also study the convergence when $h \to 0$ and we obtain a commutative diagram which relates the continuous and discrete Neumann boundary mixed elliptic optimal control problems $(P_{h\alpha}), (P_\alpha), (P_h)$ and $(P)$ by taking the limits $h \to 0$ and $\alpha \to +\infty$ respectively.

I. **Introduction**

The goal of this work is to do the numerical analysis of the convergence of the continuous Neumann boundary mixed optimal control problems with respect to a parameter (the heat transfer coeffcient) given in [15]. For distributed optimal control problems we can see [14].

We consider a bounded domain $\Omega \subset \mathbb{R}^n$ whose regular boundary $\Gamma = \partial\Omega = \Gamma_1 \cup \Gamma_2$ consists of the union of two disjoint portions $\Gamma_1$ and $\Gamma_2$ with $\text{meas}(\Gamma_1) > 0$ and $\text{meas}(\Gamma_2) > 0$. We consider the following family of continuous Neumann boundary optimal control problems $(P_\alpha)$ for each parameter $\alpha > 0$ where the control variable is the heat flux $q$ on $\Gamma_2$, that is: For each $\alpha > 0$, find the continuous Neumann boundary optimal control $q_{\alpha_{op}} \in Q = L^2(\Gamma_2)$ such that:

$$\text{Problem } (P_\alpha): \quad J_\alpha(q_{\alpha_{op}}) = \min_{q \in Q} J_\alpha(q) \tag{1}$$

where the quadratic cost functional $J_\alpha : Q \to \mathbb{R}_0^+$ is defined by the following expresion [2, 23, 30]:

$$J_\alpha(q) = \frac{1}{2}\|u_{\alpha q} - z_d\|_H^2 + \frac{M}{2}\|q\|_Q^2 \tag{2}$$

with $M > 0$ and $z_d \in H$ given, $u_{\alpha q} \in V$ is the state of the system defined by the elliptic variational equality [21]:

$$\begin{cases} a_\alpha(u_{\alpha q}, v) = (g, v)_H - (q, v)_Q + \alpha \int_{\Gamma_1} bv \, d\gamma, & \forall v \in V \\ u_{\alpha q} \in V \end{cases} \tag{3}$$

and its adjoint system state $p_{\alpha q} \in V$ is defined by the following elliptic variational equality:

$$\begin{cases} a_\alpha(p_{\alpha q}, v) = (u_{\alpha q} - z_d, v), & \forall v \in V \\ p_{\alpha q} \in V \end{cases} \tag{4}$$

where the bilinear, continuous, symmetric and coercive form $a_\alpha$ and $a$ are given by:

$$a_\alpha(u, v) = a(u, v) + \alpha \int_{\Gamma_1} uv \, d\gamma, \quad a(u, v) = \int_\Omega \nabla u \cdot \nabla v \, dx,$$
$$(u, v)_H = \int_\Omega uv \, dx, \quad (q, v)_Q = \int_{\Gamma_2} qv \, d\gamma \tag{5}$$

where $\lambda_\alpha = \lambda_1 \min(1, \alpha) > 0, \lambda_1 > 0$ and $\lambda > 0$ are the positive coercive constants of $a_\alpha, a_1$ and $a$, that is [21, 26]:

$$\lambda_\alpha \|v\|_V^2 \leq a_\alpha(v, v), \quad \forall v \in V, \text{ and } \quad \lambda \|v\|_V^2 \leq a(v, v), \quad \forall v \in V_0, \tag{6}$$

and the functional spaces are:

$$H = L^2(\Omega), \quad V = H^1(\Omega), \quad Q = L^2(\Gamma_2),$$
$$V_0 = \{v \in V, \ v/\Gamma_1 = 0\}, \quad K = \{v \in V, \ v/\Gamma_1 = b\} = b + V_0, \tag{7}$$

In (3), $g$ is the internal energy in $\Omega$, $b = Const.$ is the temperature of the external neighborhood on $\Gamma_1$, $q$ is the heat flux on $\Gamma_2$ and $\alpha > 0$ is the heat transfer coefficient on $\Gamma_1$. The systems (3) can represent the steady-state two-phase Stefan problem for adequate data [26, 27].

We also consider the following continuous Neumann boundary optimal control problem (P) where the control variable is the heat flux $q$ on $\Gamma_2$, that is: Find the continuous Neumann boundary optimal control $q_{op} \in Q$ such that:

$$\text{Problem (P):} \qquad J(q_{op}) = \min_{q \in Q} J(q) \qquad (8)$$

where the quadratic cost functional $J : Q \to \mathbb{R}_0^+$ is defined by the following expresion [2, 23, 30]:

$$J(q) = \frac{1}{2} \|u_q - z_d\|_H^2 + \frac{M}{2} \|q\|_Q^2 \qquad (9)$$

with $M > 0$ and $z_d \in H$ given, $u_q \in K$ is the state of the system defined by the following elliptic variational equality [21]:

$$\begin{cases} a(u_q, v) = (g, v)_H - (q, v)_Q, & \forall v \in V_0, \\ u_q \in K \end{cases} \qquad (10)$$

and its adjoint system state $p_q \in V$ is defined by the following elliptic variational equality:

$$\begin{cases} a(p_q, v) = (u_q - z_d, v)_H, & \forall v \in V_0 \\ p_q \in V_o \end{cases} \qquad (11)$$

In [15] the limit of the optimal control problem (1) when $\alpha \to \infty$ was studied and it was proven that:

$$\lim_{\alpha \to \infty} \|u_{\alpha q_{\alpha op}} - u_{q_{op}}\|_V = 0, \quad \lim_{\alpha \to \infty} \|p_{\alpha q_{\alpha op}} - p_{q_{op}}\|_V = 0, \quad \lim_{\alpha \to \infty} \|q_{\alpha op} - q_{op}\|_Q = 0. \qquad (12)$$

We can summary the conditions (12) saying that the Neumann boundary optimal control problems ($P_\alpha$) converges to the Neumann boundary optimal control problem (P) when $\alpha \to +\infty$.

Now, we consider the finite element method and a polygonal domain $\Omega \subset \mathbb{R}^n$ with a regular triangulation with Lagrange triangles of type 1, constituted by affine-equivalent finite element of class $C^0$ being $h$ the parameter of the finite element approximation which goes to zero [4,10]. Then, we discretize the elliptic variational equalities for the system states (3) and (10), the adjoint system states (4) and (11), and the cost functional (1) and (8) respectively. In general, the solution of a mixed elliptic boundary problem belongs to $H^r(\Omega)$ with $1 < r \leq \tfrac{3}{2} - \varepsilon$ ($\varepsilon > 0$) but there exist some examples which solutions belong to $H^r(\Omega)$ with $2 \leq r$ [1, 22, 25]. Note that mixed boundary conditions play an important role in various applications, e.g. heat conduction and electric potential problems [16].

The goal of this paper is to study the numerical analysis of the convergence (12) of the continuous Neumann boundary elliptic optimal control problems $(P_\alpha)$ to $(P)$ when $\alpha \to \infty$. The main result of this paper can be characterized by the following result:

*We have the following commutative diagram which relates the continuous and discrete Neumann boundary mixed optimal control problems $(P_{h\alpha}), (P_\alpha), (P_h)$ and $(P)$ by taking the limits $h \to 0$ and $\alpha \to +\infty$ as follows:*

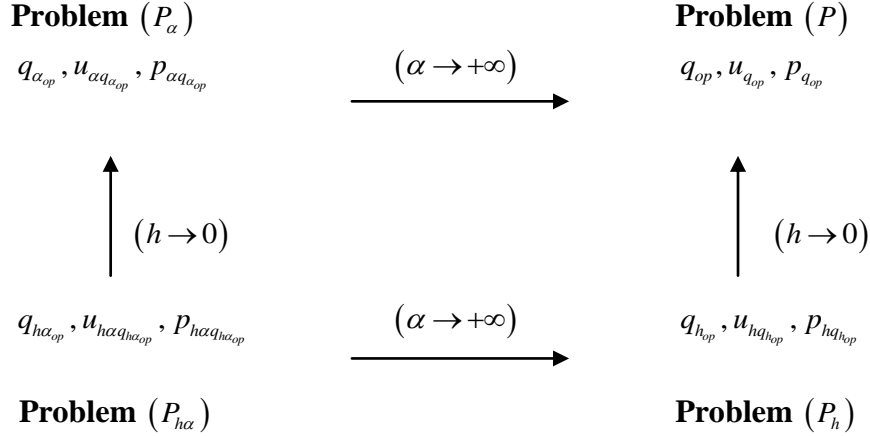

*where $u_{h\alpha q_{h\alpha op}}$ and $p_{h\alpha q_{h\alpha op}}$ are respectively the system and the adjoint system states of the discrete Neumann boundary mixed optimal control problem $(P_{h\alpha})$ for each $h > 0$ and $\alpha > 0$. Moreover, we obtain error estimates for the convergence when $h \to 0$ between the solution of problem $(P_{h\alpha})$ with respect to problem $(P_\alpha)$ for each $\alpha > 0$, and between the solution of problem $(P_h)$ with respect to problem $(P)$ respectively.*

The study of the limit $h \to 0$ of the discrete solutions of optimal control problems can be considered as a classical limit, see [3,5-9,11-13,16-20,24,28,29,31,32] but the limit $\alpha \to +\infty$, for each $h \to 0$, can be considered as a new one. Moreover, the main result given by the above commutative diagram is, from our point of view, a new and original relationship among discrete and continuous Neumann boundary mixed elliptic optimal control problems being the discrete and continuous optimal controls characterized as fixed points of certain operators.

The paper will be organized in the following manner:

In Section II we give a complement to the continuous Neumann boundary optimal control problems $(P)$ and $(P_\alpha)$ [15] by defining two contraction operators $W$ and $W_\alpha$ which allow to obtain the optimal controls $q_{op}$ and $q_{\alpha op}$ as a fixed points respctively.

In Section III we define the discrete elliptic variational equalities for the state systems $u_{hq}$ and $u_{h\alpha q}$, we define the discrete Neumann boundary cost functional $J_h$ and $J_{h\alpha}$, we define the discrete Neumann boundary optimal control problems $(P_h)$ and $(P_{h\alpha})$ and we define the discrete elliptic variational equalities for the adjoint state systems $p_{hq}$ and $p_{h\alpha q}$ for each $h > 0$ and $\alpha > 0$. We obtain properties for the optimal control problem $(P_h)$: for system state $u_{hq}$ and adjoint system state $p_{hq}$, for the discrete cost functional $J_h$ and its corresponding optimality condition. We define a contraction operator $W_h$ which allows to obtain the optimal control $q_{h op}$ as a fixed point.

We also obtain properties for the optimal control problem $(P_{h\alpha})$: for system $u_{h\alpha q}$ and adjoint system states $p_{h\alpha q}$, for the discrete cost functional $J_{h\alpha}$ and its corresponding optimality condition. We also define a contraction operator $W_{h\alpha}$ which allows to obtain the optimal control $q_{h\alpha_{op}}$ as a fixed point.

In Section IV we study the classical convergence of the discrete elliptic variational equalities for $u_{hq}, u_{h\alpha q}, p_{hq}$, and $p_{h\alpha q}$ as $h \to 0$ when $q$ is fixed (for each $\alpha > 0$). We study the convergences of the discrete optimal control problem $(P_h)$ to $(P)$ and $(P_{h\alpha})$ to $(P_\alpha)$ when $h \to 0$ (for each $\alpha > 0$). We also study the explicit error estimates for the optimal control problems $(P_h)$ and $(P_{h\alpha})$ (for each $\alpha > 0$).

In Section V we study the new convergence of the discrete Neumann boundary optimal control problems $(P_{h\alpha})$ to $(P_h)$ when $\alpha \to +\infty$ for each $h > 0$ and we obtain a commutative diagram which relates the continuous and discrete Neumann boundary mixed optimal control problems $(P_{h\alpha}), (P_\alpha), (P_h)$ and $(P)$ by taking the limits $h \to 0$ and $\alpha \to +\infty$.

In Section VI we study the convergence when $h \to 0$ of the discrete cost functional $J_h$ and $J_{h\alpha}$ corresponding to the discrete Neumann boundary optimal control problems $(P_h)$ and $(P_{h\alpha})$ respectively, $\forall \alpha > 0$.

## II. A Complement to the Continuous Neumann Boundary Optimal Control Problems $(P_\alpha)$ and $(P)$ Through Fixed Points

The unique continuous Neumann boundary optimal controls $q_{op}$ and $q_{\alpha_{op}}$ can be characterized as a fixed points on $Q$ of suitable operators $W$ and $W_\alpha$ over their optimal adjoint system states $p_{q_{op}}$ and $p_{\alpha q_{\alpha_{op}}}$ [15], for each parameter $\alpha > 0$, defined by:

$$W: Q \to Q \quad / \quad W(q) = \frac{1}{M} \gamma_0(p_q), \tag{13}$$

$$W_\alpha: Q \to Q \quad / \quad W_\alpha(q) = \frac{1}{M} \gamma_0(p_{\alpha q}), \tag{14}$$

where $\gamma_0$ is the trace operator.

**Lemma 1** *We have that:*
*(i) $W$ is a Lipschitzian operator, that is:*

$$\|W(q_2) - W(q_1)\|_Q \leq \frac{\|\gamma_0\|^2}{M \lambda^2} \|q_2 - q_1\|_Q, \quad \forall q_1, q_2 \in Q. \tag{15}$$

*(ii) $W$ is a contraction operator if and only if data M verifies the inequality*

$$M > \frac{\|\gamma_0\|^2}{\lambda^2}. \tag{16}$$

*(iii) If $M$ verifies the inequality (16) then the continuous Neumann boundary optimal control $q_{op} \in Q$ can be obtained as the unique fixed point of the operator $W$, that is:*

$$q_{op} = \frac{1}{M} \gamma_0(p_{q_{op}}) \Leftrightarrow W(q_{op}) = q_{op} \tag{17}$$

**Proof** We use the definition (13), the Lemma 3 and Corollary 5 of [15]. □

**Lemma 2** *We have that:*
*(i) $W_\alpha$ is a Lipschitzian operator, that is:*

$$\left\|W_\alpha(q_2) - W_\alpha(q_1)\right\|_Q \leq \frac{\|\gamma_0\|^2}{M\lambda_\alpha^2}\|q_2 - q_1\|_Q, \quad \forall q_1, q_2 \in Q. \tag{18}$$

*(ii) $W$ is a contraction operator if and only if data M verifies the inequality*

$$M > \frac{\|\gamma_0\|^2}{\lambda_\alpha^2}. \tag{19}$$

*(iii) If M verifies the inequality (19) then the continuous Neumann boundary optimal control $q_{\alpha_{op}} \in Q$ can be obtained as the unique fixed point of the operator $W_\alpha$, that is:*

$$q_{\alpha_{op}} = \frac{1}{M}\gamma_0\left(p_{\alpha q_{\alpha_{op}}}\right) \Leftrightarrow W_\alpha\left(q_{\alpha_{op}}\right) = q_{\alpha_{op}} \tag{20}$$

**Proof** We use the definition (14), the Lemma 8 and Corollary 10 of [15]. □

### III. Discretization by Finite Element Method and Properties

We consider the finite element method and a polygonal domain $\Omega \subset \mathbb{R}^n$ with a regular triangulation with Lagrange triangles of type 1, constituted by affine-equivalent finite element of class $C^0$ being $h$ the parameter of the finite element approximation which goes to zero [4,10]. We can take $h$ equal to the longest side of the triangles $T \in \tau_h$ and we can approximate the sets $V, V_0$ and $K$ by:

$$\begin{cases} V_h = \left\{v_h \in C^0(\overline{\Omega})/ v_h/T \in P_1(T), \forall T \in \tau_h\right\} \\ V_{0h} = \left\{v_h \in V_h / v_h/\Gamma_1 = 0\right\}; \quad K_h = b + V_{0h} \end{cases} \tag{21}$$

where $P_1$ is the set of the polymonials of degree less than or equal to 1. Let $\pi_h : V \to V_h$ be the corresponding linear interpolation operator. Then there exists a constant $c_0 > 0$ (independent of the parameter $h$) such that [4]:

$$\begin{cases} a) \|v - \pi_h(v)\|_H \leq c_0 h^r \|v\|_r, \quad \forall v \in H^r(\Omega), \quad 1 < r \leq 2 \\ b) \|v - \pi_h(v)\|_V \leq c_0 h^{r-1}\|v\|_r, \quad \forall v \in H^r(\Omega), \quad 1 < r \leq 2 \end{cases}. \tag{22}$$

We define the discrete cost functional $J_h : Q \to \mathbb{R}_0^+$ by the following expression:

$$J_h(q) = \frac{1}{2}\|u_{hq} - z_d\|_H^2 + \frac{M}{2}\|q\|_Q^2 \tag{23}$$

where $u_{hq}$ is the discrete system state defined as the solution of the following discrete elliptic variational equality [21,28,29]:

$$\begin{cases} a(u_{hq}, v_h) = (g, v_h)_H - (q, v_h)_Q, & \forall v_h \in V_{0h} \\ u_{hq} \in K_h \end{cases} \quad (24)$$

and its corresponding discrete adjoint state $p_{hq}$ is defined as the solution of the following discrete elliptic variational equality:

$$\begin{cases} a(p_{hq}, v_h) = (u_{hq} - z_d, v_h)_H, & \forall v_h \in V_{0h} \\ p_{hq} \in V_{0h}. \end{cases} \quad (25)$$

We define $u_{h0}$ as the solution of the discrete elliptic variational equality (24) for the particular case $q = 0$.

The corresponding discrete Neumann boundary optimal control problem consists in finding $q_{h_{op}} \in Q$ such that:

$$\text{Problem } (P_h): \quad J_h(q_{h_{op}}) = \underset{q \in Q}{\text{Min}} \, J_h(q). \quad (26)$$

**Lemma 3**

*(i) There exist unique solutions $u_{hq} \in K_h$ and $p_{hq} \in V_{0h}$ of the elliptic variational equalities (24), and (25) respectively $\forall g \in H, \forall q \in Q, b > 0$ on $\Gamma_1$.*

*(ii) The operator $q \in Q \to u_{hq} \in V$ is Lipschitzian, i.e.*

$$\|u_{hq_2} - u_{hq_1}\|_V \leq \frac{\|\gamma_0\|}{\lambda} \|q_2 - q_1\|_Q, \quad \forall q_1, q_2 \in Q, \forall h > 0. \quad (27)$$

*(iii) The operator $q \in Q \to p_{hq} \in V_{0h}$ is Lipschitzian and strictly monotone, i.e.*

$$-(\gamma_0(p_{hq_2}) - \gamma_0(p_{hq_1}), q_2 - q_1)_Q = \|u_{hq_2} - u_{hq_1}\|_H^2 \geq 0, \quad \forall q_1, q_2 \in Q, \forall h > 0 \quad (28)$$

$$\|p_{hq_2} - p_{hq_1}\|_V \leq \frac{1}{\lambda} \|u_{hq_2} - u_{hq_1}\|_V \leq \frac{\|\gamma_0\|}{\lambda^2} \|q_2 - q_1\|_Q, \quad \forall q_1, q_2 \in Q, \forall h > 0. \quad (29)$$

**Proof.** We use the Lax-Milgram theorem, the variational equalities (24) and (25), the coerciveness (6) and following [15,23]. □

**Theorem 4**

*(i) The discrete cost functional $J_h$ is a Q- elliptic and strictly convexe application, that is:*

$$(1-t)J_h(q_2) + tJ_h(q_1) - J_h(tq_1 + (1-t)q_2) = \frac{t(1-t)}{2} \|u_{hq_2} - u_{hq_1}\|_H^2 + M \frac{t(1-t)}{2} \|q_2 - q_1\|_Q^2$$

$$\geq M \frac{t(1-t)}{2} \|q_2 - q_1\|_Q^2, \quad \forall q_1, q_2 \in Q, \forall t \in [0,1]. \quad (30)$$

*(ii) There exists a unique optimal control $q_{h_{op}} \in Q$ that satisfies the optimization problem (26).*

*(iii) $J_h$ is a Gâteaux differenciable application and its derivative $J_h'$ is given by the following expression:*

$$J'_h(q) = Mq - \gamma_0(p_{hq}), \quad \forall q \in Q, \quad \forall h > 0. \tag{31}$$

*(iv) The optimality condition for the problem (26) is given by:*

$$J'_h(q_{h_{op}}) = 0 \Leftrightarrow q_{h_{op}} = -\frac{1}{M}\gamma_0(p_{hq_{h_{op}}}). \tag{32}$$

*(v) The operator $J'_h$ is a Lipschitzian and strictly monotone one, i.e.*

$$\|J'_h(q_2) - J'_h(q_1)\|_H \leq \left(M + \frac{\|\gamma_0\|^2}{\lambda^2}\right) \|q_2 - q_1\|_Q, \quad \forall q_1, q_2 \in Q, \forall h > 0 \tag{33}$$

$$\langle J'_h(q_2) - J'_h(q_1), q_2 - q_1 \rangle = \|u_{hq_2} - u_{hq_1}\|_H^2 + M\|q_2 - q_1\|_Q^2$$

$$\geq M\|q_2 - q_1\|_Q^2, \quad \forall q_1, q_2 \in Q, \forall h > 0 \tag{34}$$

**Proof.** We use the definition (23), the elliptic variational equalities (24) and (25) and the coerciveness (6) following [15,23]. The discrete cost functional (23) can be written as:

$$J_h(q) = \frac{1}{2} G_h(q, q) - L_h(q) + \frac{1}{2}\|u_{h0} - z_d\|_H^2, \quad \forall q \in Q \tag{35}$$

and the functional $J'_h$ is given by:

$$\langle J'_h(q), f \rangle = \lim_{t \to 0^+} \frac{J_h(q + tf) - J_h(q)}{t} = G_h(q, f) - L_h(f), \quad \forall q, f \in Q, \tag{36}$$

where the operators $G_h : Q \times Q \to \Box$, $C_h : Q \to V_{0h}$ and $L_h : Q \to \Box$ are defined by:

$$G_h(q, f) = (C_h(q), C_h(f))_H + M(q, f)_Q, \quad C_h(q) = u_{hq} - u_{h0} \tag{37}$$

$$L_h(q) = (C_h(q), z_d - u_{h0})_H \tag{38}$$

and satisfy the following property:

$$a(p_{hq}, C_h(f)) = (u_{hq} - z_d, C_h(f))_H = -(f, \gamma_0(p_{hq}))_Q, \quad \forall q, f \in Q. \quad \Box \tag{39}$$

We define the operator

$$W_h : Q \to Q \quad / \quad W_h(q) = \frac{1}{M}\gamma_0(p_{hq}). \tag{40}$$

**Theorem 5** *We have that:*
*(i) $W_h$ is a Lipschitzian operator, that is:*

$$\left\|W_h(q_2) - W_h(q_1)\right\|_Q \leq \frac{\|\gamma_0\|^2}{M\lambda^2}\|q_2 - q_1\|_Q, \quad \forall q_1, q_2 \in Q, \quad \forall h > 0. \tag{41}$$

*(ii)* $W_h$ *is a contraction operator if and only if* $M$ *is large, that is:*

$$M > \frac{\|\gamma_0\|^2}{\lambda^2}. \tag{42}$$

*(iii) If* $M$ *verifies the inequality (42) then the discrete Neumann boundary optimal control* $q_{h_{op}} \in Q$ *can be also obtained as the unique fixed point of the operator* $W_h$, *that is:*

$$q_{h_{op}} = \frac{1}{M} p_{hq_{h_{op}}} \iff W_h(q_{h_{op}}) = q_{h_{op}} \tag{43}$$

**Proof.** We use the definition (40) and the properties (29) and (32). □

We define the discrete cost functional $J_{h\alpha} : Q \to \mathbb{R}_0^+$ by the following expression:

$$J_{h\alpha}(q) = \frac{1}{2}\|u_{h\alpha q} - z_d\|_H^2 + \frac{M}{2}\|q\|_Q^2 \tag{44}$$

where $u_{h\alpha q}$ is the discrete system state defined as the solution of the following discrete elliptic variational equality [21,28,29]:

$$\begin{cases} a_\alpha(u_{h\alpha q}, v_h) = (g, v_h)_H - (q, v_h)_Q + \alpha \int_{\Gamma_1} bv_h d\gamma, & \forall v_h \in V_h \\ u_{h\alpha q} \in V_h \end{cases} \tag{45}$$

and its corresponding discrete adjoint system state $p_{h\alpha q}$ is defined as the solution of the following discrete elliptic variational equality:

$$\begin{cases} a_\alpha(p_{h\alpha q}, v_h) = (u_{h\alpha q} - z_d, v_h), & \forall v_h \in V_h \\ p_{h\alpha q} \in V_h \end{cases}. \tag{46}$$

The corresponding discrete Neumann boundary optimal control problem consists in finding $q_{h\alpha_{op}} \in Q$ such that:

$$\text{Problem } (P_{h\alpha}): \quad J_{h\alpha}(q_{h\alpha_{op}}) = \underset{q \in Q}{\text{Min}}\, J_{h\alpha}(q). \tag{47}$$

**Lemma 6**
*(i) There exist unique solutions* $u_{h\alpha q} \in V_h$ *and* $p_{h\alpha q} \in V_h$ *of the elliptic variational equalities (45), and (46) respectively* $\forall g \in H, \forall q \in Q, b > 0$ *on* $\Gamma_1$.
*(ii) The operator* $q \in Q \to u_{h\alpha q} \in V$ *is Lipschitzian, i.e.*

$$\left\|u_{h\alpha q_2} - u_{h\alpha q_1}\right\|_v \leq \frac{\|\gamma_0\|}{\lambda_\alpha}\|q_2 - q_1\|_Q, \quad \forall q_1, q_2 \in Q, \forall h > 0. \tag{48}$$

*(iii) The operator $q \in Q \to p_{h\alpha q} \in V_h$ is Lipschitzian and strictly monotone, i.e.*

$$-(p_{h\alpha q_2} - p_{h\alpha q_1}, q_2 - q_1)_Q = \|u_{h\alpha q_2} - u_{h\alpha q_1}\|_H^2 \geq 0, \quad \forall q_1, q_2 \in Q, \forall h > 0, \quad (49)$$

$$\|p_{h\alpha q_2} - p_{h\alpha q_1}\|_V \leq \frac{1}{\lambda_\alpha} \|u_{h\alpha q_2} - u_{h\alpha q_1}\|_V \leq \frac{\|\gamma_0\|}{\lambda_\alpha^2} \|q_2 - q_1\|_Q, \quad \forall q_1, q_2 \in Q, \forall h > 0. \quad (50)$$

**Proof.** We use the Lax-Milgram theorem, the variational equalities (45) and (46), the coerciveness (6) and following [15,23]. □

**Theorem 7**

*(i) The discrete cost functional $J_{h\alpha}$ is a $Q$- elliptic and strictly convexe application, that is:*

$$(1-t)J_{h\alpha}(q_2) + tJ_{h\alpha}(q_1) - J_{h\alpha}(tq_1 + (1-t)q_2) = \frac{t(1-t)}{2}\|u_{h\alpha q_2} - u_{h\alpha q_1}\|_H^2 + M\frac{t(1-t)}{2}\|q_2 - q_1\|_Q^2$$

$$\geq M\frac{t(1-t)}{2}\|q_2 - q_1\|_Q^2, \quad \forall q_1, q_2 \in Q, \forall t \in [0,1]. \quad (51)$$

*(ii) There exists a unique optimal control $q_{h\alpha_{op}} \in Q$ that satisfies the optimization problem (47).*

*(iii) $J_{h\alpha}$ is a Gâteaux differenciable aplication and its derivative $J'_{h\alpha}$ is given by the following expression:*

$$J'_{h\alpha}(q) = Mq - \gamma_0(p_{h\alpha q}), \quad \forall q \in Q, \quad \forall h > 0. \quad (52)$$

*(iv) The optimality condition for the problem (47) is given by:*

$$J'_{h\alpha}(q_{h\alpha_{op}}) = 0 \Leftrightarrow q_{h\alpha_{op}} = \frac{1}{M}\gamma_0(p_{h\alpha q_{h\alpha_{op}}}). \quad (53)$$

*(v) The application $J'_{h\alpha}$ is a Lipschitzian and strictly monotone one, i.e.*

$$\|J'_{h\alpha}(q_2) - J'_{h\alpha}(q_1)\|_Q \leq \left(M + \frac{\|\gamma_0\|^2}{\lambda_\alpha^2}\right)\|q_2 - q_1\|_Q, \quad \forall q_1, q_2 \in Q, \forall h > 0 \quad (54)$$

$$\langle J'_{h\alpha}(q_2) - J'_{h\alpha}(q_1), q_2 - q_1 \rangle = \|u_{h\alpha q_2} - u_{h\alpha q_1}\|_H^2 + M\|q_2 - q_1\|_Q^2$$

$$\geq M\|q_2 - q_1\|_Q^2, \quad \forall q_1, q_2 \in Q, \forall h > 0. \quad (55)$$

**Proof.** Similarly to the Theorem 4 we use the definition (44), the elliptic variational equalities (45) and (46) and the coerciveness (6) following [15,23]. The discrete cost functional (44) can be written as:

$$J_{h\alpha}(q) = \frac{1}{2}G_{h\alpha}(q,q) - L_{h\alpha}(q) + \frac{1}{2}\|u_{h\alpha 0} - z_d\|_H^2, \quad \forall q \in Q. \quad (56)$$

and the functional $J'_{h\alpha}$ is given by:

$$\langle J'_{h\alpha}(q), f \rangle = \lim_{t \to 0^+} \frac{J_{h\alpha}(q+tf) - J_{h\alpha}(q)}{t} = G_{h\alpha}(q,f) - L_{h\alpha}(f), \quad \forall q, f \in Q, \tag{57}$$

where the operators $G_{h\alpha} : Q \times Q \to \square$, $C_{h\alpha} : Q \to V_h$ and $L_{h\alpha} : Q \to \square$ are defined by:

$$G_{h\alpha}(q,f) = \left(C_{h\alpha}(q), C_{h\alpha}(f)\right)_H + M(q,f)_Q, \quad C_{h\alpha}(q) = u_{h\alpha q} - u_{h\alpha 0} \tag{58}$$

$$L_{h\alpha}(q) = \left(C_{h\alpha}(q), z_d - u_{h\alpha 0}\right)_H \tag{59}$$

and satisfy the following property:

$$a_\alpha\left(p_{h\alpha q}, C_{h\alpha}(f)\right) = \left(u_{h\alpha q} - z_d, C_{h\alpha}(f)\right)_H = -\left(f, \gamma_0(p_{h\alpha q})\right)_Q, \quad \forall q, f \in Q. \ \square \tag{60}$$

We define the operator

$$W_{h\alpha} : Q \to Q \quad / \quad W_{h\alpha}(q) = \frac{1}{M} \gamma_0(p_{h\alpha q}). \tag{61}$$

**Theorem 8** *We have that:*
*(i) The operator $W_{h\alpha}$ is Lipschitzian, that is:*

$$\left\| W_{h\alpha}(q_2) - W_{h\alpha}(q_1) \right\|_V \leq \frac{\|\gamma_0\|^2}{M \lambda_\alpha^2} \|q_2 - q_1\|_Q, \quad \forall q_1, q_2 \in Q, \quad \forall h > 0. \tag{62}$$

*(ii) The operator $W_{h\alpha}$ is a contraction if and only if $M$ is large, that is:*

$$M > \frac{\|\gamma_0\|^2}{\lambda_\alpha^2}. \tag{63}$$

*(iii) If $M$ verifies the inequality (63) then the discrete Neumann boundary optimal control $q_{h\alpha_{op}} \in Q$ can be also obtained as the unique fixed point of the operator $W_{h\alpha}$, that is:*

$$q_{h\alpha_{op}} = \frac{1}{M} \gamma_0(p_{h\alpha q_{h\alpha_{op}}}) \Leftrightarrow W_{h\alpha}\left(q_{h\alpha_{op}}\right) = q_{h\alpha_{op}} \tag{64}$$

**Proof.** We use the definition (61) and the properties (50) and (53). $\square$

## IV. Convergence of the Discrete Optimal Control Problems $(P_{h\alpha})$ and $(P_h)$ when $h \to 0$

We can divide the study $h \to 0$ in two parts.

### IV.1 Relationship between Neumann boundary optimal control problems $(P_h)$ and $(P)$

We obtain the following error estimations between the continuous and discrete solutions:

**Lemma 9** $\forall q \in Q$ *(fixed) we have the following properties:*

$$a(u_q - u_{hq}, v_h) = 0, \quad \forall v_h \in V_{0h} \tag{65}$$

$$a(u_q - u_{hq}, u_q - u_{hq}) \leq a(u_q - v_h, u_q - v_h), \quad \forall v_h \in K_h \tag{66}$$

$$\|u_q - u_{hq}\|_V \leq \frac{1}{\lambda} \underset{v_h \in K_h}{\text{Inf}} \|u_q - v_h\|_V . \tag{67}$$

*(ii) If the continuous system state has the regularity $u_q \in H^r(\Omega)$ $(1 < r \leq 2)$ then we have:*

$$\|u_q - u_{hq}\|_V \leq \frac{c_0}{\sqrt{\lambda}} \|u_q\|_r h^{r-1}, \quad \forall q \in Q, h > 0 . \tag{68}$$

*(iii) We have the following convergence:*

$$\lim_{h \to 0^+} \|u_q - u_{hq}\|_V = 0, \quad \forall q \in Q . \tag{69}$$

**Proof.** We use the variational equalities (10) and (24), $v_h = \pi_h(u_g)$ in the variational equality (24), the coerciveness (6) and the estimations (22). □

**Lemma 10** $\forall q \in Q$ *(fixed) we have the following properties:*

*(i)* $a(p_q - p_{hq}, \pi_h(p_q) - p_{hq}) = (u_q - u_{hq}, \pi_h(p_q) - p_{hq})$. \tag{70}

*(ii) If the continuous system state and the adjoint system state have the regularities $u_q \in H^r(\Omega)$, $p_q \in H^r(\Omega)$ $(1 < r \leq 2)$ then we have the estimations:*

$$\|p_q - p_{hq}\|_V^2 \leq c_1 \|p_q - p_{hq}\|_V h^{r-1} + c_2 h^{2r-1} \tag{71}$$

with

$$c_1 = \frac{c_0}{\lambda}\left[\|p_q\|_r + \frac{\|u_q\|_r}{\sqrt{\lambda}}\right], \quad c_2 = \frac{c_0^2}{\lambda^{3/2}} \|u_q\|_r \|p_q\|_r \tag{72}$$

and

$$\|p_q - p_{hq}\|_V \leq c_3 h^{r-1}, \quad \forall h \leq 1 \tag{73}$$

with

$$c_3 = \sqrt{c_1^2 + 2c_2} . \tag{74}$$

*(iii) We have the convergence:*

$$\lim_{h \to 0^+} \|p_q - p_{hq}\|_V = 0, \quad \forall q \in Q . \tag{75}$$

**Proof.** We use the variational equalities (11) and (25), $v_h = \pi_h(p_g)$ in the variational equality (25), the coerciveness (6) and the estimations (22). □

**Theorem 11**
*(i) If the continuous system state and adjoint system state have the regularities*

$u_{q_{op}} \in H^r(\Omega)$, $p_{q_{op}} \in H^r(\Omega)$ $(1 < r \leq 2)$ then we have the following limits:

$$\lim_{h \to 0^+} \|q_{h_{op}} - q_{op}\|_Q = 0, \quad \lim_{h \to 0^+} \|u_{hq_{h_{op}}} - u_{q_{op}}\|_V = 0, \quad \lim_{h \to 0^+} \|p_{hq_{h_{op}}} - p_{q_{op}}\|_V = 0. \tag{76}$$

**Proof.** We can divide the proof in the following steps (note that C's are positive constants independents of $h$):

(i) By using the variational equality (24) for $q = 0$ we get:

$$\|u_{h0} - b\|_V \leq \frac{1}{\lambda} \|g\|_H, \quad \forall h > 0 \text{ and } \|u_{h0}\|_V \leq C, \tag{77}$$

and therefore by using the definition of the cost functional (23) we obtain:

$$\frac{1}{2} \|u_{hq_{h_{op}}} - z_d\|_H^2 + \frac{M}{2} \|q_{h_{op}}\|_Q^2 \leq \frac{1}{2} \|u_{h0} - z_d\|_H^2 \leq C,$$

that is

$$\|u_{hq_{h_{op}}}\|_H \leq C, \quad \|q_{h_{op}}\|_Q \leq C, \quad \forall h > 0. \tag{78}$$

(ii) By using the variational equality (24) we have:

$$\|u_{hq_{h_{op}}} - b\|_V \leq \frac{1}{\lambda}\left[\|g\|_H + \|q\|_Q \|\gamma_0\|\right] \leq C, \quad \forall h > 0,$$

and then

$$\|u_{hq_{h_{op}}}\|_V \leq C, \quad \forall h > 0. \tag{79}$$

(iii) By using the variational equality (25) we have:

$$\|p_{hq_{h_{op}}}\|_V \leq \frac{1}{\lambda} \|u_{hq_{h_{op}}} - z_d\|_H \leq C, \quad \forall h > 0. \tag{80}$$

(iv) From the above estimations we have that:

$$\begin{cases} a) & \exists f \in Q / q_{h_{op}} \to f \text{ in Q weak as } h \to 0^+ \\ b) & \exists \eta \in V / u_{hq_{h_{op}}} \to \eta \text{ in } V \text{ weak (in } H \text{ strong) as } h \to 0^+ \\ c) & \exists \xi \in V / p_{hq_{h_{op}}} \to \xi \text{ in } V \text{ weak (in } H \text{ strong) as } h \to 0^+ \end{cases} \tag{81}$$

(v) By using the above three weak convergences we can pass to the limit as $h \to 0^+$, and we obtain by uniqueness of the variational equalities (24) and (25) that: $\eta = u_f$, $\xi = p_f$ and $f = q_{op}$.

(vi) By the other hand, by using (6) and the variational equality (24) we have:

$$\lambda \|u_{hq_{h_{op}}} - u_{q_{op}}\|_V^2 \leq \left(g, u_{q_{op}} - u_{hq_{h_{op}}}\right)_H + \left(q_{h_{op}} - q_{op}, u_{hq_{h_{op}}} - b\right)_Q - (q_{op}, u_{q_{op}} - u_{hq_{h_{op}}})_Q \to 0, \text{ as } h \to 0$$

and therefore we deduce that:

$$\lim_{h \to 0^+} \left\| u_{hq_{h_{op}}} - u_{q_{op}} \right\|_V = 0. \tag{82}$$

By using (6) and the variational equality (25) we have:

$$\lambda \left\| p_{hq_{h_{op}}} - p_{q_{op}} \right\|_V^2 \leq \left( u_{hq_{h_{op}}} - u_{q_{op}}, p_{hq_{h_{op}}} \right)_H - a\left( p_{q_{op}}, p_{hq_{h_{op}}} - p_{q_{op}} \right) \to 0, \text{ as } h \to 0^+,$$

and then we deduce that:

$$\lim_{h \to 0^+} \left\| p_{hq_{h_{op}}} - p_{q_{op}} \right\|_V = 0. \tag{83}$$

(vii) By using the definition (23) we can pass to the limit as $h \to 0^+$ and we deduce that:

$$\lim_{h \to 0^+} \left\| q_{h_{op}} \right\|_Q = \left\| q_{op} \right\|_Q. \tag{84}$$

(viii) From the weak convergence (90) and the property (84) we deduce that:

$$\lim_{h \to 0} \left\| q_{h_{op}} - q_{op} \right\|_Q = 0, \tag{85}$$

and all the limits (76) hold. □

**Remark 1.** If M verifies the inequality (42) we can obtain that $f = g_{op}$ by using the characterization of the fixed point (43), and then we obtain $f = \frac{1}{M} p_f$ when $h \to 0$. By uniqueness of the optimal control $q_{op} \in Q$ we deduce that $f = q_{op}$.

**Theorem 12** *If* M *verifies the inequality (42) and the continuous system state and adjoint system state have the regularities $u_{q_{op}} \in H^r(\Omega)$, $p_{q_{op}} \in H^r(\Omega)(1 < r \leq 2)$ then we have the following error bonds:*

$$\left\| q_{h_{op}} - q_{op} \right\|_Q \leq Ch^{r-1}, \left\| u_{hq_{h_{op}}} - u_{q_{op}} \right\|_V \leq Ch^{r-1}, \left\| p_{hq_{h_{op}}} - p_{q_{op}} \right\|_V \leq Ch^{r-1} \tag{86}$$

*where C's are constants independents of h.*

**Proof.** By using the fixed point property (43) we have:

$$\left\| q_{h_{op}} - q_{op} \right\|_Q \leq \frac{\|\gamma_0\|}{M} \left[ \left\| p_{q_{op}} - p_{hq_{op}} \right\|_V + \left\| p_{hq_{op}} - p_{hq_{h_{op}}} \right\|_V \right] \leq \frac{\|\gamma_0\|}{M} \left[ c_3 h^{r-1} + \frac{1}{\lambda^2} \left\| q_{h_{op}} - q_{op} \right\|_Q \right],$$

that is

$$\left\| q_{h_{op}} - q_{op} \right\|_Q \leq \frac{\lambda^2 c_3}{\|\gamma_0\| \left( \frac{M\lambda^2}{\|\gamma_0\|} - 1 \right)} h^{r-1}, \quad \forall h \in (0,1). \tag{87}$$

By using the variational equalities (10) and (24), we have:

$$a\left(u_{hq_{h_{op}}} - u_{q_{op}}, v_h\right) = -\left(q_{h_{op}} - q_{op}, v_h\right)_Q, \quad \forall v_h \in V_{0h}. \tag{88}$$

Therefore, by using (6) and (88), we get:

$$\lambda \left\|u_{hq_{h_{op}}} - u_{q_{op}}\right\|_V^2 \leq a(u_{hq_{h_{op}}} - u_{q_{op}}, u_{hq_{h_{op}}} - u_{q_{op}}) = a(u_{hq_{h_{op}}} - u_{q_{op}}, u_{hq_{h_{op}}} - \pi_h(u_{q_{op}}) + \pi_h(u_{q_{op}}) - u_{q_{op}})$$

$$\leq \left\|u_{hq_{h_{op}}} - u_{q_{op}}\right\|_V \left(\left\|q_{h_{op}} - q_{op}\right\|_Q + \left\|\pi_h(u_{q_{op}}) - u_{q_{op}}\right\|_V\right) + \left\|q_{h_{op}} - q_{op}\right\|_Q \left\|\pi_h(u_{q_{op}}) - u_{q_{op}}\right\|_Q$$

$$\leq \left\|u_{hq_{h_{op}}} - u_{q_{op}}\right\|_V \left(\frac{\lambda^2 c_3}{\frac{M\lambda^2}{\|\gamma_0\|^2} - 1} h^{r-1} + c_0 \left\|u_{q_{op}}\right\|_r h^{r-1}\right) + \frac{\lambda^2 c_3}{\frac{M\lambda^2}{\|\gamma_0\|^2} - 1} h^{r-1} c_0 \left\|u_{q_{op}}\right\|_r h^r$$

$$= \left\|u_{hq_{h_{op}}} - u_{q_{op}}\right\|_V \lambda c_4 h^{r-1} + \lambda c_5 h^{2r-1}$$

that is

$$\left\|u_{hq_{h_{op}}} - u_{q_{op}}\right\|_V^2 \leq c_4 \left\|u_{hq_{h_{op}}} - u_{q_{op}}\right\|_V h^{r-1} + c_5 h^{2r-1} \tag{89}$$

where

$$c_4 = \frac{c_3 \lambda}{\frac{M\lambda^2}{\|\gamma_0\|^2} - 1} + \frac{c_0}{\lambda} \left\|u_{q_{op}}\right\|_r, \quad c_5 = \frac{c_3 \lambda}{\frac{M\lambda^2}{\|\gamma_0\|^2} - 1} c_0 \left\|u_{q_{op}}\right\|_r.$$

Therefore from the above inequality (89) we deduce that

$$\left\|u_{hq_{h_{op}}} - u_{q_{op}}\right\|_V \leq c_6 h^{r-1}, \quad \forall h \leq 1, \text{ with } c_6 = \sqrt{c_4^2 + 2c_5}. \tag{90}$$

By using the variational equalities (11) and (25), we have:

$$a\left(p_{hq_{h_{op}}} - p_{q_{op}}, v_h\right) = \left(u_{hq_{h_{op}}} - u_{q_{op}}, v_h\right)_H, \quad \forall v_h \in V_{0h}. \tag{91}$$

If we take $v_h = \pi_h\left(p_{q_{op}}\right) - p_{hq_{h_{op}}} \in V_{0h}$ in (91), in a similar way to the previous result, we can deduce:

$$\left\|p_{hq_{h_{op}}} - p_{q_{op}}\right\|_V^2 \leq c_7 \left\|p_{hq_{h_{op}}} - p_{q_{op}}\right\|_V h^{r-1} + c_8 h^{2r-1}, \quad \forall h \leq 1 \tag{92}$$

with the constants

$$c_7 = \frac{c_6 + c_0 \left\|p_{q_{op}}\right\|_r}{\lambda}, \quad c_8 = \frac{c_0 c_6}{\lambda} \left\|p_{q_{op}}\right\|_r,$$

and therefore we obtain the inequality

$$\left\| p_{hq_{h_{op}}} - p_{q_{op}} \right\|_V \leq c_9 h^{r-1}, \quad \forall h \leq 1, \text{ with } c_9 = \sqrt{c_7^2 + 2c_8}, \tag{93}$$

and the thesis holds. □

## IV.2 Relationship between Neumann boundary optimal control problems $(P_{h\alpha})$ and $(P_\alpha)$

Following the above section we can obtain the following error estimations between the continuous and discrete solutions of the Neumann boundary optimal control problems $(P_{h\alpha})$ and $(P_\alpha)$.

**Lemma 13** *(i) If the continuous system state and adjoint system state have the regularities $u_{\alpha q} \in H^r(\Omega)$, $p_{\alpha q} \in H^r(\Omega)$ $(1 < r \leq 2)$ then $\forall \alpha > 0, \forall q \in Q$ we have the estimations:*

$$\left\| u_{h\alpha q} - u_{\alpha q} \right\|_V \leq ch^{r-1}, \quad \left\| p_{h\alpha q} - p_{\alpha q} \right\|_V \leq ch^{r-1} \tag{94}$$

*where the constants $c$'s are independents of h.*
*(ii) We have the following limits:*

$$\lim_{h \to 0^+} \left\| u_{h\alpha q} - u_{\alpha q} \right\|_V = 0, \quad \lim_{h \to 0^+} \left\| p_{h\alpha q} - p_{\alpha q} \right\|_V = 0, \quad \forall \alpha > 0, \forall q \in Q. \tag{95}$$

**Proof**. In a similar way to the one developed in Lemmas 9 and 10 and by using the variational equalities (3), (4), (45) and (46), the thesis holds. □

**Theorem 14**
*(i) If the continuous system state and adjoint system state have the regularities $u_{\alpha q_{\alpha_{op}}}, p_{\alpha q_{\alpha_{op}}} \in H^r(\Omega)$ $(1 < r \leq 2)$ and the inequality $\dfrac{M\lambda_1^2}{\|\gamma_0\|^2} > 1$ is verified then we have the following estimations $\forall \alpha > 1, \forall q \in Q$:*

$$\left\| q_{h\alpha_{op}} - q_{\alpha_{op}} \right\|_Q \leq ch^{r-1}, \quad \left\| u_{h\alpha q_{h\alpha_{op}}} - u_{\alpha q_{\alpha_{op}}} \right\|_V \leq ch^{r-1}, \quad \left\| p_{h\alpha q_{h\alpha_{op}}} - p_{\alpha q_{\alpha_{op}}} \right\|_V \leq ch^{r-1} \tag{96}$$

*where the constants $c$'s are independents of h.*
*(ii) We have the following limits:*

$$\lim_{h \to 0^+} \left\| q_{h\alpha_{op}} - q_{\alpha_{op}} \right\|_Q = 0, \quad \lim_{h \to 0^+} \left\| u_{h\alpha q_{h\alpha_{op}}} - u_{\alpha q_{\alpha_{op}}} \right\|_V = 0, \lim_{h \to 0^+} \left\| p_{h\alpha q_{h\alpha_{op}}} - p_{\alpha q_{\alpha_{op}}} \right\|_V = 0, \quad \forall \alpha > 1 \tag{97}$$

**Proof**. In a similar way to the one developed in Theorems 11 and 12, and by using the variational equalities (3), (4), (45) and (46), the thesis holds. □

**Remark 2**. The restriction $\alpha > 1$ can be replaced by $\alpha_0 \leq \alpha$ for any $\alpha_0 > 0$.

# V. Convergence of the Discrete Optimal Control Problems $(P_{h\alpha})$ when $\alpha \to +\infty$

For a fixed $h > 0$ we have:

**Lemma 15** For a fixed $q \in Q$ we have the following limits:

$$\lim_{\alpha \to +\infty} \|u_{h\alpha q} - u_{\alpha q}\|_V = 0, \forall q \in Q, \forall h > 0, \qquad (98)$$

$$\lim_{\alpha \to +\infty} \|p_{h\alpha q} - p_{hq}\|_V = 0, \forall q \in Q, \forall h > 0. \qquad (99)$$

**Proof.** For fixed $q \in Q, h > 0$, and by using the variational equalities (3) and (45), and by splitting the bilinear form $a_\alpha$, when $\alpha > 1$, by [26,29]

$$a_\alpha(u,v) = a_1(u,v) + (\alpha - 1) \int_{\Gamma_1} uv \, d\gamma \quad , \qquad (100)$$

we obtain the following estimations:

$$\|u_{h\alpha q} - u_{hq}\|_V \leq c, \quad (\alpha - 1) \int_{\Gamma_1} (u_{h\alpha q} - b)^2 d\gamma \leq c, \quad \forall \alpha > 1. \qquad (101)$$

From the above inequalities (101) we deduce that:

$$\exists \eta_{hq} \in V / \begin{cases} u_{h\alpha q} \longrightarrow \eta_{hq} & \text{in } V \text{ weak (in } H \text{ strong) as } \alpha \to +\infty \\ \eta_{hq}/\Gamma_1 = b \end{cases}. \qquad (102)$$

By using the variational equality (45) we can pass to the limit when $\alpha \to +\infty$, and by uniqueness of the variational equality (24) we obtain that $\eta_{hq} = u_{hq}$. By using the above properties, and the variational equalities (3) y (45), we deduce that:

$$u_{h\alpha q} \to u_{hq} \text{ in } V \text{ strong as } \alpha \to +\infty. \qquad (103)$$

Finally, by using a similar method developed before for the discrete system state we can obtain the limit $\alpha \to +\infty$ for the discrete adjoint system state, i.e. (99) holds. □

**Theorem 16** We have the following limits:

$$\lim_{\alpha \to +\infty} \|u_{h\alpha q_{h\alpha op}} - u_{hq_{h op}}\|_V = 0, \quad \forall h > 0, \qquad (104)$$

$$\lim_{\alpha \to +\infty} \|p_{h\alpha q_{h\alpha op}} - p_{hq_{h op}}\|_V = 0, \quad \forall h > 0, \qquad (105)$$

$$\lim_{\alpha \to +\infty} \|q_{h\alpha op} - q_{h op}\|_Q = 0, \quad \forall h > 0. \qquad (106)$$

**Proof.** From now on we consider a fixed parameter $h > 0$ and we also consider that c's represent positive constants independents of $\alpha > 0$. If we use the variational equality (45) for the particular case $q = 0$ and we splitting the bilinear form (100) we obtain the following estimations:

$$\|u_{h\alpha 0} - u_{h0}\|_V \leq c, \quad (\alpha - 1) \int_{\Gamma_1} (u_{h\alpha 0} - b)^2 d\gamma \leq c, \quad \forall \alpha > 1. \qquad (107)$$

From the definition of the discrete optimal control problem (44) we obtain the following estimations:

$$\frac{1}{2}\left\|u_{h\alpha q_{h\alpha_{op}}} - z_d\right\|_H^2 + \frac{M}{2}\left\|q_{h\alpha_{op}}\right\|_Q^2 \leq \frac{1}{2}\left\|u_{h\alpha 0} - z_d\right\|_H^2 \leq c, \quad \forall \alpha > 0,$$

and therefore we deduce the estimations:

$$\left\|u_{h\alpha q_{h\alpha_{op}}}\right\|_H \leq c, \quad \left\|q_{h\alpha_{op}}\right\|_Q \leq c, \quad \forall \alpha > 0. \tag{108}$$

Now, by using the variational equality (45) for the optimal state system and splitting the bilinear form (100) we get the estimations:

$$\left\|u_{h\alpha q_{h\alpha_{op}}} - u_{hq_{h_{op}}}\right\|_V \leq c, \quad (\alpha-1)\int_{\Gamma_1} \left(u_{h\alpha q_{h\alpha_{op}}} - b\right)^2 d\gamma \leq c, \quad \forall \alpha > 1. \tag{109}$$

In a similar way by using the variational equality (46) for the discrete adjoint state system we deduce the following estimations:

$$\left\|p_{h\alpha q_{h\alpha_{op}}} - p_{hq_{h_{op}}}\right\|_V \leq c, \quad (\alpha-1)\int_{\Gamma_1} p_{h\alpha q_{h\alpha_{op}}}^2 d\gamma \leq c, \quad \forall \alpha > 1. \tag{110}$$

Then, from the above properties we have that:

$$\exists\, f_h \in H \,/\, q_{h\alpha_{op}} \longrightarrow f_h \text{ in } Q \text{ weak as } \alpha \to +\infty \tag{111}$$

$$\exists\, \eta_h \in V \,/\, \begin{cases} u_{h\alpha q_{h\alpha_{op}}} \longrightarrow \eta_h & \text{in } V \text{ weak (in } H \text{ strong) as } \alpha \to +\infty \\ \eta_h / \Gamma_1 = b \end{cases} \tag{112}$$

$$\exists\, \xi_h \in V \,/\, \begin{cases} p_{h\alpha q_{h\alpha_{op}}} \longrightarrow \xi_h & \text{in } V \text{ weak (in } H \text{ strong) as } \alpha \to +\infty \\ \xi_h / \Gamma_1 = 0 \end{cases} \tag{113}$$

By using the three above weak convergences we can pass to the limit $\alpha \to +\infty$, and by uniqueness of the variational equalities (24) and (25) we get that $\eta_h = u_{hf_h}$, $\xi_h = p_{hf_h}$. By using (23) and (44) we can pass to the limit $\alpha \to +\infty$, and by uniqueness of the discrete optimal control problem (26) we have $f_h = q_{h_{op}}$. Therefore, we deduce that

$$\eta_h = u_{hf_h} = u_{hq_{h_{op}}}, \quad \xi_h = p_{hf_h} = p_{hq_{h_{op}}}. \tag{114}$$

By using the variational equalities (3) and (45) for the discrete system state, and the variational equalities (4) and (46) for the discrete adjoint system state, we obtain the following strong convergences:

$$\lim_{\alpha \to +\infty}\left\|u_{h\alpha q_{h\alpha_{op}}} - u_{hq_{h_{op}}}\right\|_V = 0, \quad \lim_{\alpha \to +\infty}\int_{\Gamma_1}\left(u_{h\alpha q_{h\alpha_{op}}} - b\right)^2 d\gamma = 0, \quad \forall h > 0, \tag{115}$$

and

$$\lim_{\alpha \to +\infty} \left\| p_{h\alpha q_{h\alpha op}} - p_{hq_{hop}} \right\|_V = 0, \quad \lim_{\alpha \to +\infty} \int_{\Gamma_1} p_{h\alpha q_{h\alpha op}}^2 \, d\gamma = 0, \quad \forall h > 0. \tag{116}$$

On the other hand, we can pass to the limit $\alpha \to +\infty$ in the discrete cost functional (23) and (44), and we obtain:

$$\lim_{\alpha \to +\infty} \left\| q_{h\alpha_{op}} \right\|_Q = \left\| q_{h_{op}} \right\|_Q, \quad \forall h > 0. \tag{117}$$

From this result (117) and the weak convergence of the discrete optimal controls we obtain the strong convergence of the optimal control, that is:

$$\lim_{\alpha \to +\infty} \left\| q_{h\alpha_{op}} - q_{h_{op}} \right\|_Q = 0, \quad \forall h > 0. \qquad \square \tag{118}$$

## VI. Convergence of the Discrete Cost Functional when $h \to 0$

Following Section IV.1 we have:

**Lemma 17** *If* M *verifies the inequality (42) and the continuous system state and adjoint system state have the regularities* $u_{q_{op}} \in H^r(\Omega), p_{q_{op}} \in H^r(\Omega) (1 < r \leq 2)$ *then we have the following error bonds:*

$$\frac{M}{2} \left\| q_{h_{op}} - q_{op} \right\|_Q^2 \leq J(q_{h_{op}}) - J(q_{op}) \leq Ch^{2(r-1)} \tag{119}$$

$$\frac{M}{2} \left\| q_{h_{op}} - q_{op} \right\|_H^2 \leq J_h(q_{op}) - J_h(q_{h_{op}}) \leq Ch^{2(r-1)} \tag{120}$$

$$\left| J_h(q_{op}) - J(q_{op}) \right| \leq Ch^{r-1} \tag{121}$$

$$\left| J_h(q_{h_{op}}) - J(q_{op}) \right| \leq Ch^{r-1} \tag{122}$$

*where C's are constants independents of h.*

**Proof.** Estimations (119) and (120) follow from the estimations (66) and (96), and the equalities:

$$J(q_{h_{op}}) - J(q_{op}) = \frac{1}{2} \left\| u_{q_{h_{op}}} - u_{q_{op}} \right\|_H^2 + \frac{M}{2} \left\| q_{h_{op}} - q_{op} \right\|_Q^2, \tag{123}$$

$$J_h(q_{op}) - J_h(q_{h_{op}}) = \frac{1}{2} \left\| u_{hq_{h_{op}}} - u_{hq_{op}} \right\|_H^2 + \frac{M}{2} \left\| q_{h_{op}} - q_{op} \right\|_Q^2. \tag{124}$$

Estimation (121) follows from the estimations (27), (66) and (86), and the inequality:

$$\left| J_h(q) - J(q) \right| \leq \left( \frac{1}{2} \left\| u_{hq} - u_q \right\|_H + \left\| u_q - z_d \right\|_H \right) \left\| u_{hq} - u_q \right\|_H, \quad \forall q \in Q. \tag{125}$$

Finally, estimation (122) follows from the previous results and the triangular inequality for norms.

□

**Remark 3** We can also obtain for the optimal control problem $(P_{h\alpha})$ similar results to the one given in Lemma 17, e.g.

$$\left| J_{h\alpha}(q_{op}) - J_\alpha(q_{op}) \right| \leq Ch^{r-1}, \tag{126}$$

$$\left| J_{h\alpha}(q_{h\alpha_{op}}) - J_\alpha(q_{\alpha_{op}}) \right| \leq Ch^{r-1}, \tag{127}$$

which proof will be omitted here.

## Conclusions

We have studied the numerical analysis of the discrete Neumann boundary optimal control problems $(P_h)$ and $(P_{h\alpha})$, and the corresponding asymptotic behaviour when $\alpha \to \infty$ and $h \to 0$ by using the finite element method. We have defined the discrete cost functional $J_h$ and $J_{h\alpha}$, the discrete variational equalities for the system states $u_{hg}$ and $u_{h\alpha g}$ for each $\alpha, h > 0$, and the discrete variational equalities for the adjoint system states $p_{hg}$ and $p_{h\alpha g}$ for each $\alpha, h > 0$. We have characterized the discrete Neumann boundary optimal control heat fluxes $q_{h_{op}}$ and $q_{h\alpha_{op}}$ as a fixed point on $Q$ of suitable discrete operators $W_h$ and $W_{h\alpha}$ over his adjoint system states $p_{hg_{op}}$ and $p_{h\alpha g_{h\alpha_{op}}}$ respectively for each $\alpha > 0$. We have also studied the convergence of the discrete Neumann boundary optimal control problems $(P_{h\alpha})$ to $(P_h)$ when $\alpha \to \infty$ for each $h > 0$, and when $h \to 0$ for each $\alpha > 0$, and we have obtained a commutative diagram (see Introduction) which relates the continuous and discrete Neumann boundary mixed optimal control problems $(P_{h\alpha}), (P_\alpha), (P_h)$ and $(P)$ by taking the limits $h \to 0$ and $\alpha \to \infty$.

## Acknowledgements
This paper has been partially sponsored by the Project PIP No. 534 from CONICET – Univ. Austral (Rosario, Argentina), and Grant AFOSR-SOARD FA9550-14-1-0122.